**Hankel determinants of backward shifts of the coefficients of a partial theta function.**

Johann Cigler

**Abstract.**

We study some polynomials which are related to Hankel determinants of backward shifts of the coefficients of a partial theta function. In this version an appendix is added which gives a simple formula for the coefficients of the reciprocal of the partial theta function.

## 1. Introduction

Consider the double sequence $\left(a(n,q)\right)_{n\in\mathbb{Z}}$ with $a(n,q)=q^{\binom{n}{2}}$ for $n\geq 0$ and $a(n,q)=0$ for $n<0$ and the Hankel determinants

(1) $$D_{-m,n}(q)=\det\left(a(-m+i+j,q)\right)_{i,j=0}^{n-1}$$

for $m\leq 0.$ For $0<n\leq m$ we get $D_{-m,n}(q)=0$ because the first row of the matrix $\left(a(-m+i+j,q)\right)_{i,j=0}^{n-1}$ vanishes. For $n=0$ we set $D_{-m,0}(q)=1$ by definition.

Since $D_{0,n+1}(q)$ does not vanish we can write

(2) $$D_{-m,n+m+1}(q)=(-1)^{\binom{m+1}{2}}r_{m,n}(q)q^{m\binom{n}{2}}D_{0,n+1}(q)$$

with a uniquely determined function $r_{m,n}(q).$

Note that the generating function for the sequence $\left(a(n,q)\right)$ is the partial theta function $\sum_{n\geq 0}q^{\binom{n}{2}}x^n.$

Computations suggest the

**Conjecture**

*The functions* $r_{m,n}(q)$ *are monic polynomials with integer coefficients with* $\deg r_{m,n}(q)=\dfrac{mn(n+m+2)}{2}$ *which satisfy* $r_{m,n}(1)=1$ *and* $r_{m,n}(0)=(-1)^{mn}.$

For example,

$$\left(r_{m,n}(q)\right)_{m,n=0}^{2}=\begin{pmatrix}1 & 1 & 1\\ 1 & -1+q+q^2 & 1-q-q^2+q^4+q^5\\ 1 & 1-2q-q^2+q^3+q^4+q^5 & 1-2q-q^2+2q^3+2q^4+2q^5-3q^6-2q^7-2q^8+2q^{10}+q^{11}+q^{12}\end{pmatrix}$$

We shall prove this conjecture for some special cases and derive some recurrences for the general case.

Let us first consider the Hankel determinants

(3) $$d_{m,n}(q)=\det\left(q^{\binom{m+i+j}{2}}\right)_{i,j=0}^{n-1}$$

for $m\in\mathbb{Z}$ of the double sequence $\left(q^{\binom{n}{2}}\right)_{n\in\mathbb{Z}}$. We get

$$d_{m,n}(q)=\det\left(q^{\binom{m+i+j}{2}}\right)_{i,j=0}^{n-1}=\det\left(q^{\binom{i}{2}+\binom{j}{2}+\binom{m}{2}+im+jm+ij}\right)_{i,j=0}^{n-1}=q^{n\binom{m}{2}}\prod_{i=0}^{n-1}q^{\binom{i}{2}}\prod_{j=0}^{n-1}q^{\binom{j}{2}}\prod_{i=0}^{n-1}q^{im}\prod_{j=0}^{n-1}q^{jm}\det\left(q^{ij}\right)_{i,j=0}^{n-1}$$
$$=q^{2\binom{n}{3}+2m\binom{n}{2}+n\binom{m}{2}}\det\left(q^{ij}\right)_{i,j=0}^{n-1}.$$

The Vandermonde determinant evaluation (cf. [3], (2.1)) $\det\left(x_i^j\right)_{i,j=0}^{n-1}=\prod_{0\le i<j\le n-1}\left(x_j-x_i\right)$ gives

$$\det\left(q^{ij}\right)_{i,j=0}^{n-1}=\det\left(\left(q^i\right)^j\right)_{i,j=0}^{n-1}=\prod_{0\le i<j\le n-1}\left(q^j-q^i\right)=q^{\binom{n}{3}}\prod_{0\le i<j\le n-1}\left(q^{j-i}-1\right)=q^{\binom{n}{3}}(q-1)^{\binom{n}{2}}\prod_{j=0}^{n-1}[j]!.$$

Therefore, we get

(4) $$d_{m,n}(q)=q^{3\binom{n}{3}+2m\binom{n}{2}+n\binom{m}{2}}(q-1)^{\binom{n}{2}}\prod_{j=0}^{n-1}[j]!$$

with the usual $q-$ notation $[n]=[n]_q=1+q+\cdots+q^{n-1}=\dfrac{q^n-1}{q-1}$ and $[n]!=[n]_q!=[1][2]\cdots[n].$

**2. The polynomials $r_{m,n}(q)$ for small $m$.**

**Theorem 1**

(5) $$r_{1,n}(q)=q^{\binom{n+2}{2}}-(q-1)^{n+1}[n+1]!=\sum_{k=0}^{n}(q-1)^k[k]!q^{\binom{n+2}{2}-\binom{k+2}{2}}.$$

**Proof**

To compute $D_{-1,n+2}(q)$ we first use the expansion of a determinant by minors: $\det\left(a_{i,j}\right)_{i,j=0}^{n-1}=\sum_{j=0}^{n-1}(-1)^j a_{0,j}\det A_{0,j}$ where the minors $A_{0,j}$ are obtained by crossing out the first row and $j-$th column.

Since $a(-1,q)=q^{\binom{-1}{2}}=q$ we get

(6) $$D_{-1,n+2}(q)=d_{-1,n+2}(q)-qd_{1,n+1}(q).$$

For example

$$D_{-1,4}(q)=\det\begin{pmatrix}0&1&1&q\\1&1&q&q^3\\1&q&q^3&q^6\\q&q^3&q^6&q^{10}\end{pmatrix}=\det\begin{pmatrix}q&1&1&q\\1&1&q&q^3\\1&q&q^3&q^6\\q&q^3&q^6&q^{10}\end{pmatrix}-q\det\begin{pmatrix}1&q&q^3\\q&q^3&q^6\\q^3&q^6&q^{10}\end{pmatrix}.$$

By (4) we get $\dfrac{d_{-1,n+2}(q)}{d_{0,n+1}(q)}=\dfrac{q^{3\binom{n+2}{3}-2\binom{n+2}{2}+(n+2)\binom{-!}{2}}(q-1)^{\binom{n+2}{2}}\prod_{j=0}^{n+1}[j]!.}{q^{3\binom{n+1}{3}}(q-1)^{\binom{n+1}{2}}\prod_{j=0}^{n}[j]!.}=q^{\binom{n}{2}}(q-1)^{n+1}[n+1]!$

and

$$\frac{qd_{1,n+1}(q)}{d_{0,n+1}(q)}=\frac{q^{3\binom{n+1}{3}+2\binom{n+1}{2}+1}(q-1)^{\binom{n+1}{2}}\prod_{j=0}^{n}[j]!}{q^{3\binom{n+1}{3}}(q-1)^{\binom{n+1}{2}}\prod_{j=0}^{n}[j]!}=q^{\binom{n+2}{2}+\binom{n}{2}}$$

Thus $\dfrac{D_{-1,n+2}(q)}{d_{0,n+1}(q)}=q^{\binom{n}{2}}(q-1)^{n+1}[n+1]!-q^{\binom{n+2}{2}+\binom{n}{2}}$ or

(7) $$r_{1,n}(q)=-\frac{D_{-1,n+2}(q)}{d_{0,n+1}(q)}\frac{1}{q^{\binom{n}{2}}}=q^{\binom{n+2}{2}}-(q-1)^{n+1}[n+1]!.$$

Another way to compute $r_{1,n}(q)$ uses Dodgson's condensation theorem (cf. [3], Prop. 10) which gives

(8) $$D_{m,n+2}(q)D_{m+2,n}(q)-D_{m+2,n+1}(q)D_{m,n+1}(q)+D_{m+1,n+1}(q)^2=0.$$

Setting $g(n)=D_{-1,n}(q)$ gives

$g(n+2)=\dfrac{-D_{0,n+1}(q)^2+D_{1,n+1}(q)g(n+1)}{D_{1,n}(q)}$ and thus

$$r_{1,n}(q)=\frac{d_{0,n+1}(q)^2+d_{1,n+1}(q)r_{1,n-1}(q)q^{\binom{n-1}{2}}d_{0,n}(q)}{d_{1,n}(q)q^{\binom{n}{2}}d_{0,n+1}(q)}=(q-1)^n[n]!+q^{n+1}r_{1,n-1}(q).$$

which gives

(9) $$r_{1,n}(q)=\sum_{k=0}^{n}(q-1)^k[k]!q^{\binom{n+2}{2}-\binom{k+2}{2}}.$$

Formulae (7) and (9) show that $r_{1,n}(q)$ satisfies the recurrence

(10) $$r_{1,n}(q)=q^{n+1}r_{1,n-1}(q)+(q-1)^n[n]!$$

with $r_{1,0}(q)=1.$

By induction we see that $\deg\left(r_{1,n}(q)\right)=\frac{n(n+3)}{2}.$ Note that

$\deg\left((q-1)^n[n]!\right)=\frac{n(n+1)}{2}<n+1+\frac{(n-1)(n+2)}{2}=\frac{n(n+3)}{2}.$

For $q=2$ we get

$\left(r_{1,n}(2)\right)_{n\geq 0}=\left(1,5,43,709,23003,1481957,190305691,48796386661,\cdots\right),$

which occurs in OEIS [4], A114604, in another context.

**Theorem 2**

*For $m=2$ we get by condensation*

(11) $$r_{2,n}(q)=r_{1,n}^2(q)+\left(q^{n+1}-1\right)q^{n+2}r_{2,n-1}(q).$$

or equivalently

(12) $$r_{2,n}(q)=f(n,q)\sum_{j=0}^{n}\frac{r_{1,j}^2(q)}{f(j,q)}$$

with $f(n,q)=(q-1)^{n+1}[n+1]!q^{\frac{(n+1)(n+4)}{2}}.$

By induction we see that $\deg\left(r_{2,n}(q)\right)=n(n+4).$

For the proof let $h(n)=D_{-2,n}(q).$ Then we get $h(n+3)=\frac{-\left(D_{-1,n+2}(q)\right)^2+D_{0,n+2}(q)h(n+2)}{D_{0,n+1}(q)}$

and thus

$$r_{2,n}(q)=-\frac{h(n+3)}{q^{2\binom{n}{2}}D_{0,n+1}(q)}=\frac{-\left(D_{-1,n+2}(q)\right)^2-r_{2,n-1}(q)q^{2\binom{n-1}{2}}D_{0,n+2}(q)D_{0,n}(q)}{q^{2\binom{n}{2}}D_{0,n+1}(q)^2}$$

$$=-\frac{r_{1,n}^2(q)q^{2\binom{n}{2}}d_{0,n+1}(q)^2+r_{2,n-1}(q)q^{2\binom{n-1}{2}}d_{0,n+2}(q)D_{0,n}(q)}{q^{2\binom{n}{2}}d_{0,n+1}(q)^2}=r_{1,n}^2(q)+\left(q^{n+1}-1\right)q^{n+2}r_{2,n-1}(q).$$

## 3. The polynomials $r_{m,n}(q)$ for small $n$.

Consider the matrices $V_{k,n}(q)=\left(a(k-n+i+j,q)\right)_{i,j=0}^{n-1}$ with $v_{k,n}(q)=\det V_{k,n}(q)$. Note that in $V_{k,n}(q)$ there are $k$ non-vanishing entries in the first row.

For example

$$V_{1,4}(q)=\begin{pmatrix}0&0&0&1\\0&0&1&1\\0&1&1&q\\1&1&q&q^3\end{pmatrix},\quad V_{2,4}(q)=\begin{pmatrix}0&0&1&1\\0&1&1&q\\1&1&q&q^3\\1&q&q^3&q^6\end{pmatrix}.$$

It is clear that $v_{1,n}(q)=(-1)^{\binom{n}{2}}$.

From

$$v_{k,n+k}(q)=D_{-n,n+k}(q)=D_{-n,n+1+(k-1)}(q)=(-1)^{\binom{n+1}{2}}r_{n,k-1}(q)q^{n\binom{k-1}{2}}D_{0,k}(q)$$
$$=(-1)^{\binom{n+1}{2}}q^{(n+k)\binom{k-1}{2}}r_{n,k-1}(q)(q-1)^{\binom{k}{2}}\prod_{j=0}^{k-1}[j]!$$

we get $v_{2,m+2}(q)=(-1)^{\binom{m+1}{2}}r_{m,1}(q)(q-1)$ and $v_{3,m+3}(q)=(-1)^{\binom{m+1}{2}}q^{(m+3)}r_{m,2}(q)(q-1)^3(q+1)$.

Let us first compute the polynomials $r_{m,1}(q)$.

**Theorem 3**

*Let* $\sum_{n\geq 0}u(n,q)x^n=\dfrac{1}{\sum_{n\geq 0}q^{\binom{n}{2}}x^n}$. *Then*

(13) $$r_{m,1}(q)=\frac{u(m+2,q)}{1-q}.$$

**Proof.**

For $n \ge 2$ $V_{2,n}(q)$ is obtained from $V_{1,n+1}(q)$ by deleting the first row and column. By Cramer's rule

$\left(V_{1,n+1}(q)\right)^{-1} = \dfrac{1}{\det\left(V_{1,n+1}(q)\right)}\left(\alpha_{j,i}\right)_{i,j=0}^{n}$ with $\alpha_{j,i} = (-1)^{i+j}\det A_{j,i}$, where $A_{i,j}$ is the matrix obtained by crossing out row $i$ and column $j$ in $V_{1,n+1}(q)$. Thus $A_{0,0} = V_{2,n}(q)$.

Therefore $(-1)^{\binom{n+1}{2}} v_{2,n}(q)$ is the entry in position $(0,0)$ of the inverse matrix of $V_{1,n+1}(q)$.

It is easy to verify that

$$\left(V_{1,n+1}(q)\right)^{-1} = \left(u(n-i-j,q)\right)_{i,j=0}^{n}. \tag{14}$$

For example

$$\left(V_{1,4}\right)^{-1} = \begin{pmatrix} u(3,q) & u(2,q) & u(1,q) & u(0,q) \\ u(2,q) & u(1,q) & u(0,q) & 0 \\ u(1,q) & u(0,q) & 0 & 0 \\ u(0,q) & 0 & 0 & 0 \end{pmatrix} = \begin{pmatrix} -1+2q-q^3 & 1-q & -1 & 1 \\ 1-q & -1 & 1 & 0 \\ -1 & 1 & 0 & 0 \\ 1 & 0 & 0 & 0 \end{pmatrix}.$$

Therefore

$$v_{2,n}(q) = (-1)^{\binom{n+1}{2}} u(n,q). \tag{15}$$

The first terms of $u(n,q)$ are

$$\left(u(n,q)\right)_{n\ge 0} = \left(1,-1,1-q,-1+2q-q^3,1-3q+q^2+2q^3-q^6,-1+4q-3q^2-3q^3+2q^4+2q^6-q^{10},\cdots\right).$$

Since $v_{2,n+2}(q) = (-1)^{\binom{n+1}{2}}(q-1)r_{n,1}(q)$ we get (13).

The polynomials $u(n,q) \in \mathbb{Z}[q]$ have degree $\binom{n}{2}$. This follows from

$\sum_{j=0}^{n} u(n-j,q)q^{\binom{j}{2}} = 0$ for $n > 0$ by induction.

Therefore, $\deg r_{n,1}(q) = \deg u(n+2,q) - 1 = \dfrac{n(n+3)}{2}$.

$$\left(r_{n,1}(q)\right)_{n\ge 0} = \left(1,-1+q+q^2,1-2q-q^2+q^3+q^4+q^5,-1+3q-3q^3-q^4-q^5+q^6+q^7+q^8+q^9,\cdots\right).$$

**Remark**

A simple formula for $u(n,q)$ will be given in the Appendix.

To compute $v_{k,n+k}(q)$ we can use the following

**Lemma** (cf. [1] Theorem 2, [2] Prop. 2.5):

*Let* $s(x)=\sum_{n\geq 0} s_n x^n$ *with* $s_0=1$ *and* $t(x)=\frac{1}{s(x)}=\sum_{n\geq 0} t_n x^n$.

*Setting* $s_n=t_n=0$ *for* $n<0$ *we get for* $M\in\mathbb{N}$

(16) $$\det\left(s_{i+j-M}\right)_{i,j=0}^{N+M}=(-1)^{N+\binom{M+1}{2}}\det\left(t_{i+j+M+2}\right)_{i,j=0}^{N-1}.$$

Choosing $s(x)=\sum_{n\geq 0} q^{\binom{n}{2}}x^n$ we get

(17) $$v_{k,n+k}(q)=\det\left(a(-n+i+j,q)\right)_{i,j=0}^{n-1}=(-1)^{k-1+\binom{n+1}{2}}\det\left(u(i+j+n+2,q)\right)_{i,j=0}^{k-1}.$$

For $k=1$ and $k=2$ this gives again $v_{1,n+1}(q)=(-1)^{\binom{n}{2}}$ and $v_{2,n+2}(q)=(-1)^{\binom{n-1}{2}}u(n+2,q)$.

Using these special cases we get by condensation

(18) $$v_{k,n+k}(q)v_{k,n+k-2}(q)-v_{k-1,n+k-1}(q)v_{k+1,n+k-1}(q)+v_{k,n+k-1}^2(q)=0.$$

For $n\geq 2$ this implies

(19) $$r_{m,n}(q)=\frac{1}{\left(q^n-1\right)q^{m+n+1}r_{m+2,n-2}(q)}\det\begin{pmatrix} r_{m,n-1}(q) & r_{m+1,n-1}(q)\\ r_{m+1,n-1}(q) & r_{m+2,n-1}(q)\end{pmatrix}.$$

## Appendix

After the first version had been posted I found a simple formula for $u(n,q)$. I want to thank Michael Schlosser for valuable hints.

From $f(x)=\sum_{n\geq 0} q^{\binom{n}{2}}x^n=1+xf(qx)$ we get

(20) $$\sum_{n\geq 0} u(n,q)x^n=\frac{1}{f(x)}=\frac{1}{1+xf(qx)}=\sum_{k\geq 0}(-1)^k x^k f(qx)^k.$$

Define $q-$analogs $\left\langle {n \atop k} \right\rangle_q$ of the binomial coefficients by the recursion

(21) $$\left\langle {n \atop k} \right\rangle_q=\sum_{j=0}^{\min(n,k)} q^{\binom{j+2}{2}}\left\langle {n-1-j \atop k-j} \right\rangle_q$$

with initial values $\left\langle {n \atop 0} \right\rangle_q = q^{n+1}$ for $n \ge -1$ and $\left\langle {-1 \atop k} \right\rangle_q = 0$ for $k > 0.$

The corresponding $q-$Pascal triangle begins with

$$\left(\left\langle {n \atop k} \right\rangle_q\right)_{n,k=0}^{5} = \begin{pmatrix} q & 0 & 0 & 0 & 0 & 0 \\ q^2 & q^3 & 0 & 0 & 0 & 0 \\ q^3 & 2q^4 & q^6 & 0 & 0 & 0 \\ q^4 & 3q^5 & q^6+2q^7 & q^{10} & 0 & 0 \\ q^5 & 4q^6 & 3q^7+3q^8 & 2q^9+2q^{11} & q^{15} & 0 \\ q^6 & 5q^7 & 6q^8+4q^9 & q^9+6q^{10}+3q^{12} & q^{12}+2q^{13}+2q^{16} & q^{21} \end{pmatrix}.$$

**Remark**

The sequence of coefficients $1,1,1,1,2,1,1,3,1,2,1,\cdots$ occurs in OEIS, A260533.

By induction we get

$$x^k f(qx)^k = \sum_{n \ge k} \left\langle {n-1 \atop n-k} \right\rangle_q \left(\frac{x}{q}\right)^n \tag{22}$$

since

$$x^k f(qx)^k = xf(qx)x^{k-1}f(qx)^{k-1} = \sum_{i \ge 0} q^{\binom{i+1}{2}} x^{i+1} \sum_{j \ge 0} \left\langle {j+k-1 \atop j} \right\rangle_q \left(\frac{x}{q}\right)^{j+k}$$

$$= \sum_{i,j} x^{i+j+k+1} q^{\binom{i+1}{2}-j-k} \left\langle {j+k-1 \atop j} \right\rangle_q = \sum_{m \ge 0} x^{m+k+1} \sum_{i+j=m} q^{\binom{i+2}{2}-i-1-j-k} \left\langle {j+k-1 \atop j} \right\rangle_q$$

$$= \sum_{m \ge 0} \frac{x^{m+k+1}}{q^{m+k+1}} \sum_{i=0}^{m} q^{\binom{i+2}{2}} \left\langle {m-i+k-1 \atop m-i} \right\rangle_q = \sum_{n \ge k+1} \left(\frac{x}{q}\right)^n \sum_{i=0}^{n-k-1} q^{\binom{i+2}{2}} \left\langle {n-i-2 \atop n-k-1-i} \right\rangle_q = \sum_{n \ge k+1} \left(\frac{x}{q}\right)^n \left\langle {n-1 \atop n-k-1} \right\rangle_q.$$

By (20) we finally get

$$u(n,q) = \frac{1}{q^n} \sum_{k=0}^{n-1} (-1)^{n-k} \left\langle {n-1 \atop k} \right\rangle_q. \tag{23}$$

Thus, the polynomials $u(n,q)$ are essentially alternating sums of the entries of the rows of the $q-$Pascal triangle $\left(\left\langle {n \atop k} \right\rangle_q\right)$. For example,

$$u(6,q) = 1 - 5q + 6q^2 + 3q^3 - 6q^4 - 2q^6 + 2q^7 + 2q^{10} - q^{15}$$
$$= \frac{1}{q^6}\left(q^6 - 5q^7 + 6q^8 + 4q^9 - (q^9 + 6q^{10} + 3q^{12}) + q^{12} + 2q^{13} + 2q^{16} - q^{21}\right).$$

**Remark**

Michael Schlosser [5] conjectured the following combinatorial interpretation of $\left\langle {n \atop k} \right\rangle_q$.

**Conjecture**

*For an integer partition* $\lambda = \left(\lambda_1, \lambda_2, \cdots, \lambda_\ell\right)$ *of* $n = \lambda_1 + \lambda_2 + \cdots + \lambda_\ell$ *with* $\lambda_1 \geq \lambda_2 \geq \cdots \geq \lambda_\ell$ *let* $w(\lambda) = coef\left(\lambda\right) q^{ex(\lambda)}$ *with* $coef(\lambda) = \prod_{i=1}^{\ell-1} \binom{\lambda_i}{\lambda_{i+1}}$ *and* $ex(\lambda) = \sum_{i=1}^{\ell} i\lambda_i$.

*Denoting by* $P_{n.k}$ *the set of all partitions of* $n$ *with first term* $k$ *we get*

(24)
$$\left\langle {n \atop k} \right\rangle_q = \sum_{\lambda \in P_{n+1,n+1-k}} w(\lambda).$$

For example $\left\langle {5 \atop 3} \right\rangle_q = \sum_{\lambda \in P_{6.3}} w(\lambda) = w\left(3,3\right) + w(3,2,1) + w(3,1,1,1)$

$$= \binom{3}{3} q^3 + \binom{3}{2}\binom{2}{1} q^{3+4+3} + \binom{3}{1}\binom{1}{1}\binom{1}{1} q^{3+2+3+4} = q^3 + 6q^{10} + 3q^{12}.$$